\documentclass[11pt,a4paper]{article}

\usepackage{amsmath,amssymb}

\newcommand{\counte}{theorem}

\allowdisplaybreaks

\renewcommand{\thefootnote}{\fnsymbol{footnote}}

\begin{document}

\renewcommand{\thefootnote}{\arabic{footnote}}

\centerline{\Large\bf A proof of Toponogov's theorem}
\centerline{\Large\bf in Alexandrov geometry\footnote{Supported by NSFC 11971057 and BNSF Z190003. \hfill{$\,$}}}

\vskip5mm

\centerline{Shengqi Hu, Xiaole Su,  Yusheng Wang\footnote{The
corresponding author (E-mail: wyusheng@bnu.edu.cn). \hfill{$\,$}}}

\vskip6mm

\noindent{\bf Abstract.} This paper aims to give an elementary proof for Toponogov's theorem in Alexandrov geometry with lower curvature bound. The idea of the proof comes from the fact that, in Riemannian geometry, sectional curvature can be embodied in the second variation formula.

\vskip1mm


\vskip1mm


\vskip6mm

\setcounter{section}{-1}

\section {Introduction}

An Alexandrov space $X$ with curvature $\geq k$
is roughly defined to be a locally complete intrinsic metric space
on which Toponogov's type theorem holds locally ([BGP]).
In fact, if $X$ is complete, then such a type theorem
holds globally on $X$, which is called Toponogov's theorem in Alexandrov geometry.
In its proofs ([BGP], [AKP]), to simplify the arguments, it is always assumed that
there exists a minimal geodesic (shortest path) between any two points in $X$
\footnote{If $X$ is both complete and locally compact,
there is a shortest path between any two points in it.}.
For any two minimal geodesics $[pq]$ and $[qr]$ in $X$, we can define the angle $\angle pqr$ naturally;
and then, Toponogov's theorem in Alexandrov geometry can be formulated as follows  ([BGP]).

\vskip2mm

\noindent {\bf Theorem A} {\it Let $X$ be a complete Alexandrov space with curvature $\geq k$.
Then for any triangle $\triangle pqr \subset X$ (a union of three minimal geodesics $[pq], [pr], [qr]$), we have that
$$\angle pqr\geq\tilde \angle_k pqr.\eqno{(0.1)}$$}
\noindent\hskip5mm In Theorem A, $\tilde \angle_k pqr$ denotes the angle $\angle \tilde p\tilde q\tilde r$ in the comparison triangle
$\triangle \tilde p\tilde q\tilde r\subset\Bbb S^2_k$ of $\triangle pqr$ (i.e. $|\tilde p\tilde q|=|pq|,|\tilde p\tilde
r|=|pr|$, $|\tilde q\tilde r|=|qr|$), where $\mathbb S^2_k$ denotes the complete and simply connected $2$-dimensional space form with constant curvature $k$.
And in a proof of Theorem A, one need only to consider the case where $|pq|+|pr|+|qr|<\frac{2\pi}{\sqrt k}$ (and thus each of $|pq|, |pr|, |qr|$ is less than $\frac{\pi}{\sqrt{k}}$) if $k>0$ ([BGP])
\footnote{Once Theorem A has been proven, it can be shown that $|pq|\leq\frac{\pi}{\sqrt{k}}$ and $|pq|+|pr|+|qr|\leq\frac{2\pi}{\sqrt k}$
for all $p,q,r\in X$ ([BGP]).
Moreover, if $|pq|=\frac{\pi}{\sqrt{k}}$, it is a convention that $\tilde \angle_k pqr=0$.}.

So far, there are several proofs for Theorem A ([BGP], [Pl], [Sh], [Wa]), each of which is skillful and the one in
[Pl] is quite beautiful. In Riemannian geometry, sectional curvature can be seen just from
geodesic variations (the second variation formula).
However, all known proofs of Theorem A have no direct relation with such an idea.
The main goal of this paper is to provide a proof of such an idea for Theorem A.

Our strategy is: If $\angle pqr<\tilde\angle_k pqr$ for some $\triangle pqr\subset X$, we can split $\triangle pqr$
into two triangles along some $[ps]$ with $s\in[qr]^\circ$ (the interior part of $[qr]$)
so that $\angle psq<\tilde\angle_k psq$  or $\angle psr<\tilde\angle_k psr$.
By repeating this process, we can locate an $o\in [qr]$ such that there is information
of ``curvature $\not\geq k$'' around some $[po]$, i.e.
there is $\triangle pr_1r_2$ with $r_i$ sufficiently close to $o$ such that $\angle pr_1r_2<\tilde\angle_k pr_1r_2$.
Note that Theorem A holds on a small ball $B(o,\delta)$ (i.e. (0.1) is true for triangles in $B(o,\delta)$) because $X$ is of curvature
$\geq k$. Then via Alexandrov's Lemma (Lemma 1.3 below) on $\triangle pr_1r_2$,
we can locate another $o'\in B(o,\delta)$ with $|po'|\leq|po|-\frac{\delta}{3}$ such that
there is still information of ``curvature $\not\geq k$'' around some $[po']$.
Step by step, such information can be transmitted to a small neighbourhood of $p$, a contradiction (because Theorem A holds around $p$).

\vskip2mm

In the rest of the paper, $X$ always denotes the space in Theorem A.

\section{Tools of the proof}

Note that a proof of Theorem A has to depend only on its local version, where 
we have the following basic property on angles ([BGP]).

\vskip2mm

\noindent {\bf Lemma 1.1} {\it Let $[pq],[rr']\subset X$  with $r\in[pq]^\circ$. Then
$\angle prr'+\angle qrr'=\pi$.}

\vskip2mm

Furthermore, we have the following easy observation (cf. [AKP], [Wa]).

\vskip2mm

\noindent {\bf Lemma 1.2} {\it Let $[pq],[qr]\subset X$. Then for $q_i\in [qr]$
with $q_i\to q$ as $i\to\infty$, $$|pq_i|\leq |pq|-\cos\angle pqr\ \cdot |qq_i|+o(|qq_i|)\
\footnote{In addition, if $X$ is locally compact, then
we can select $[pq]$ such that (1.1) is an equality (the first variation formula) once Theorem A has been proven ([BGP], [AKP]).}
.\eqno{(1.1)}$$}
\noindent{\it Proof}. Since Theorem A holds around $q$,
there is $\bar p\in [pq]$ near $q$ such that
$$|\bar pq_i|\leq|\bar pq|-\cos\angle pqr\ \cdot |qq_i|+o(|qq_i|).$$
Then (1.1) follows from that $|pq_i|\leq |p\bar p|+|\bar pq_i|$.
\hfill $\Box$

\vskip2mm

Besides Lemmas 1.1 and 1.2, we will also use Alexandrov's lemma ([BGP]).

\vskip2mm

\noindent {\bf Lemma 1.3} {\it Let $\triangle pqr$, $\triangle pqs$, $\triangle abc \subset\Bbb S^2_k$
(where $\triangle pqr$ and $\triangle pqs$ are joined to each other in an exterior way along $[pq]$) such that
$|ab|=|pr|$, $|ac|=|ps|$, $|bc|=|qr|+|qs|$, and $|ab|+|ac|+|bc|<\frac{2\pi}{\sqrt k}$ if $k>0$. Then
$\angle pqr+\angle pqs\leq \pi$ (resp. $\geq\pi$) if and only if  $\angle prq\geq\angle abc$ and  $\angle psq\geq\angle acb$
(resp. $\angle prq\leq\angle abc$ and  $\angle psq\leq\angle acb$).}

\section{Proof of Theorem A}

Due to the similarity of proofs for $k<, =, >0$, we only consider the case where $k=0$. And for simpleness,
we denote by $\tilde \angle pqr$ the angle $\tilde \angle_0 pqr$.

In the proof, we will argue by contradiction, and  say that an angle $\angle pqr$ is {\it bad} if
$\angle pqr<\tilde \angle pqr$. If each angle of a triangle is not bad,  we call the triangle a {\it good} one.

First of all, we have the following observation about a `bad' angle just via Lemmas 1.1 and 1.2 (cf. [Wa], [SSW]).

\vskip2mm

\noindent {\bf Lemma 2.1}{\it\ For a triangle $\triangle pr_1r_2\subset X$,
if $\angle pr_1r_2$ is bad, then there is  $s_0\in [r_1r_2]^\circ$ such that for any $[ps_0]$
\begin{equation}\label{eqn2.1}
\text{ $\angle ps_0r_1$  or $\angle ps_0r_2$ is bad};
\end{equation}
in particular, for each $i$,
\begin{equation}\label{eqn2.2}
\text{if $|r_is_0|\leq|r_ip|$, then $\angle ps_0r_i$ is bad}.
\end{equation}
Moreover, we have that \begin{equation}\label{eqn2.3}
|ps_0|<\max_{i=1,2}\{|pr_i|\}.
\end{equation}}
\noindent{\it Proof}. Let $\triangle \tilde p\tilde r_1\tilde r_2\subset\Bbb R^2$ be the comparison triangle of $\triangle pr_1r_2$.
By Lemma 1.2 (and the first variation formula on $\Bbb R^2$), the badness of $\angle pr_1r_2$ implies that the function $|ps|-|\tilde p\tilde s|$ with $s\in[r_1r_2], \tilde s\in[\tilde r_1\tilde r_2]$ and $|r_is|=|\tilde r_i\tilde s|$
attains a negative minimum at some $s_0\in [r_1r_2]^\circ$. By Lemma 1.2 again, for any $[ps_0]$, we have that
$$
\angle ps_0r_i\geq\angle \tilde p\tilde s_0\tilde r_i,\ \ i=1,2.
$$
It then has to hold that
$$\angle ps_0r_i=\angle \tilde p\tilde s_0\tilde r_i,\ \ i=1,2	$$
because $\angle ps_0r_1+\angle ps_0r_2=\pi$ by Lemma 1.1.
On the other hand, since $|ps_0|<|\tilde p\tilde s_0|$, there is $\tilde p'\in[\tilde p\tilde s_0]^\circ$
such that $|\tilde p'\tilde s_0|=|ps_0|$. It is clear that $|\tilde p'\tilde r_i|<|\tilde p\tilde r_i|$ for at least one of $i$, which implies $\angle \tilde p\tilde s_0\tilde r_i<\tilde \angle ps_0r_i$, i.e. (2.1) holds. Especially, it is easy to see that $|\tilde p'\tilde r_i|<|\tilde p\tilde r_i|$ if $|r_is_0|\leq|r_ip|$\ \footnote{For the case where $k>0$ and $|r_ip|\geq\frac{\pi}{2\sqrt{k}}$, it needs $\angle\tilde r_i\tilde p\tilde s_0<\frac\pi2$ besides $|r_is_0|\leq|r_ip|$ (note that $\triangle \tilde p\tilde r_1\tilde r_2\subset\Bbb S_k^2$ and $|r_i p|<\frac{\pi}{\sqrt{k}}$, cf. the comments following Theorem A).
In fact, `$\angle\tilde r_i\tilde p\tilde s_0<\frac\pi2$' holds obviously if $|r_is_0|\ll|r_ip|$ (here, it needs only `$|r_is_0|\leq|r_ip|$'
for other cases including $k\leq 0$).}; and thus (2.2) follows.
Moreover, it is clear that $|\tilde p\tilde s_0|<\max\limits_{i=1,2}\{|pr_i|\}$ \footnote{For the case where $k>0$ and $\max\limits_{i=1,2}\{|pr_i|\}\geq\frac{\pi}{2\sqrt{k}}$, it needs to be modified to `$|\tilde p\tilde s_0|<\max\limits_{i=1,2}\{|pr_i|\}+\tau(|r_1r_2|)$',
where  $\tau(|r_1r_2|)/|r_1r_2|\to 0$ as $|r_1r_2|\to0$; moreover, the function $\tau(\cdot)$ can be chosen to be the same one (i.e. not depending on $|pr_i|$) if $\max\limits_{i=1,2}\{|pr_i|\}\leq c<\frac{\pi}{\sqrt{k}}$ for some constant $c$.}, which implies (\ref{eqn2.3}).
\hfill$\Box$

\vskip2mm
\noindent {\bf Corollary 2.2} {\it\ For $\triangle pr_1r_2\subset X$,
if $\angle pr_1r_2$ is bad, then there is $\bar s\in [r_1r_2]$ such that
\begin{equation}\label{eqn2.4}
|p\bar s|\leq \max_{i=1,2}\{|pr_i|\},
\end{equation}
and for any small $\delta>0$ there exist $s_1,s_2\in [r_1r_2]\cap B(\bar s,\delta)$ such that for any $[ps_i]$
\begin{equation}\label{eqn2.5}
\text{$\angle ps_1s_2$ or $\angle ps_2s_1$ is bad.}
\end{equation}}
\noindent{\it Proof}. By Lemma 2.1, there is $[ps_0]$ with $s_0\in [r_1r_2]^\circ$ such that
$\angle ps_0r_1$ or $\angle ps_0r_2$ is bad, say $\angle ps_0r_1$.
Then we can apply Lemma 2.1 to $\triangle ps_0r_1$ again. By repeating this infinite times,
we can locate a desired $\bar s$ if we take into account (2.2) and (2.3).
\hfill$\Box$

\vskip2mm

In our proof of Theorem A,  Lemmas 1.3 and 2.1 shall be the mere keys. For simpleness, we first consider the case
where $X$ is, in addition, locally compact.

\vskip2mm

\noindent{\it Proof of Theorem A where $X$ is locally compact}.

We argue by contradiction.  Assume that there is a $\triangle pqr$ such that
$\angle pqr<\tilde \angle pqr$. Then by Corollary 2.2 (see (2.5)), we can consider the nonempty set
$$\mathcal{S}\triangleq \{x\in X| \ \forall\ \delta>0, \ \exists\ \triangle pr_1r_2 \text{ with } r_i\in B(x,\delta) \text{ s.t. }\text{$\angle pr_1r_2$ is bad}  \}.$$
It is clear that $p\not\in\mathcal{S}$ because Theorem A holds around $p$, and that $\mathcal{S}$ is closed.
Note that a closed and bounded subset in $X$ is compact because $X$ is complete and locally compact (cf. Chapter 2 in [BBI]).
So, there is $\bar o\in \mathcal{S}$ such that $|p\bar o|=\min_{o\in\mathcal{S}}\{|po|\}>0$.

On the other hand, Theorem A holds around any $o\in \mathcal{S}$, and thus we can define a positive function $\delta(o)\triangleq\min\{\frac{|po|}{2}, \delta_{o}\}$, where $\delta_{o}$ is the maximal number such that
$$\text{any $\triangle xyz\subset B(o,\delta_o)$ is good}.\eqno{(2.6)}$$
{\bf Claim}: In $B(o,\delta(o))$, there is another point $o'\in\mathcal{S}$ such that $|po'|\leq |po|-\frac{\delta(o)}{3}$.

It is obvious that the claim contradicts the existence of $\bar o$. Thereby, we just need to verify the claim. Due to $o\in \mathcal{S}$, there is a $\triangle pr_1r_2$ with $r_i\in B(o,\delta(o))$ such that
$$|or_i|\ll \delta(o), \text{ and $\angle pr_1r_2$ is bad}.$$
Let $\bar r_1\in [pr_1]$ with $|p\bar r_1|=|po|-\frac{\delta(o)}{3}.$
Note that for any $[\bar r_1r_2]$,  $\triangle r_2r_1\bar r_1$ is good because it is contained in $B(o,\delta(o))$; so via Lemmas 1.3 and 1.1
on $\triangle pr_1r_2$, the badness of $\angle pr_1r_2$ implies
$$\text{$\angle p\bar r_1r_2$ is bad}.$$ Then
we can apply (2.2) in Lemma 2.1 (note that $|\bar r_1 r_2|<\delta(o)<|p\bar r_1|$ because $\delta(o)\leq\frac{|po|}{2}$) to $\triangle p\bar r_1r_2$ to locate an $s\in [\bar r_1r_2]^\circ$ such that
$$\text{$\angle ps\bar r_1$ is bad for any $[ps]$}.$$
Put $\bar r_2\triangleq s$ if $|ps|\leq|po|-\frac{\delta(o)}{3}$; otherwise, similar to $\bar r_1$, we can select $\bar r_2\in [ps]$ with
$$|p\bar r_2|=|po|-\frac{\delta(o)}{3}, \text{ and $\angle p\bar r_2\bar r_1$ is bad for any $[\bar r_1\bar r_2]$}.$$
(Here, $\bar r_1, s, \bar r_2$ all lie in $B(o,\epsilon)$ with $\epsilon$ very close to $\frac{\delta(o)}{3}$. Note that
$|r_2x|+|xp|$ with $x=\bar r_1, s, \bar r_2$ is very close to $|po|$ because $|or_i|\ll \delta(o)$.)
Then by Corollary 2.2, the badness of $\angle p\bar r_2\bar r_1$ enables us to locate a point $o'\in [\bar r_1\bar r_2]\cap\mathcal{S}$, which lies in
$B(o,\delta(o))$, such that $|po'|\leq |po|-\frac{\delta(o)}{3}$ (see (2.4) and note that $|p\bar r_i|\leq|po|-\frac{\delta(o)}{3}$
\footnote{For the case where $k>0$ and $\max\limits_{i=1,2}\{|p\bar r_i|\}\geq\frac{\pi}{2\sqrt{k}}$,
$\delta(o)$ should be additionally so small that $\tau(\delta(o))\ll\delta(o)$ (e.g. $\tau(\delta(o))<0.01\delta(o)$), where $\tau(\cdot)$ is just the function in Footnote 7.}).
\hfill$\Box$

\vskip2mm

Eventually, we prove Theorem A for general $X$, i.e. $X$ might not be locally compact.

\vskip2mm

\noindent{\it Proof of Theorem A where $X$ is not locally compact}.

When $X$ is not locally compact, there might be no $\bar o\in \mathcal{S}$ with
$|p\bar o|=\min_{o\in\mathcal{S}}\{|po|\}$ in the proof right above. However, starting with an $o_1\in\mathcal{S}$
and by the claim below (2.6), we can step by step obtain $\{o_i\}_{i=1}^\infty\subset\mathcal{S}$
such that
$$|po_{i+1}|\leq |po_i|-\frac{\delta(o_i)}{3} \text{ and } |o_io_{i+1}|<\delta(o_i).$$
It follows that
$$|po_{i+1}|\leq |po_1|-\sum_{j=1}^i\frac{\delta(o_j)}{3}, \text{ and thus }\sum_{j=1}^i\frac{\delta(o_j)}{3}\leq|po_1|.$$
This implies that $\lim\limits_{i\to\infty}\delta(o_i)=0$, and that $\{o_i\}_{i=1}^\infty$ is a Cauchy sequence (note that $|o_io_{i+1}|<\delta(o_i)$).
Then by the completeness of $X$ and the closedness of $\mathcal{S}$, $\{o_i\}_{i=1}^\infty$ has a limit point $\bar o\in \mathcal{S}$. However, it follows that $\lim\limits_{i\to\infty}\delta(o_i)=\delta(\bar o)>0$, a contradiction.
\hfill$\Box$


\vskip5mm

\noindent School of Mathematical Sciences (and Lab. math. Com.
Sys.), Beijing Normal University, Beijing, 100875
P.R.C.

\noindent E-mail: 13716456647@163.com; suxiaole$@$bnu.edu.cn; wyusheng$@$bnu.edu.cn

\end{document}